\documentclass[oneside,english]{amsart}
\usepackage[T1]{fontenc}
\usepackage[latin9]{inputenc}
\usepackage{amsthm}
\usepackage{amssymb}

\numberwithin{figure}{section} 
\theoremstyle{plain}
\theoremstyle{plain}
\newtheorem{thm}{Theorem}
  \theoremstyle{plain}
  \newtheorem{lem}[thm]{Lemma}
  \theoremstyle{remark}
  \newtheorem{rem}[thm]{Remark}

\usepackage{babel}

\begin{document}

\title{Interpolation of compact Lipschitz operators}

\author{Michael Cwikel, Alon Ivtsan and Eitan Tadmor}

\address{M.C.~and A.I.: Department of Mathematics, Technion - Israel Institute
of Technology, Haifa 32000, Israel. }

\address{E.T.: Department of Mathematics, University of Maryland, CSCAMM,
4149 CSIC Building \#406, Paint Branch Drive, College Park, MD 20742-2389,
USA. }

\email{mcwikel@math.technion.ac.il, aloniv@techunix.technion.ac.il \textit{and}
\vskip -0.01mm
tadmor@cscamm.umd.edu}

\thanks{The research of the first named author was supported by the Technion
V.P.R.\ Fund and by the Fund for Promotion of Research at the Technion. }

\subjclass[2000]{Primary 46B70, Secondary 47H99, 46B50}

\keywords{Nonlinear operators, Lipschitz operators, compact operators, interpolation
space. }

\maketitle
\begin{abstract}
Let $\left(A_{0},A_{1}\right)$ and $\left(B_{0},B_{1}\right)$ be
Banach couples such that $A_{0}\subset A_{1}$ and $\left(B_{0},B_{1}\right)$
satisfies Arne Persson's approximation condition (H). Let $T:A_{1}\to B_{1}$
be a possibly nonlinear Lipschitz mapping which also maps $A_{0}$
into $B_{0}$ and satisfies the quantitative compactnesss condition
$Ta\in\left\Vert a\right\Vert _{A_{0}}K$ for each $a\in A_{0}$,
where $K$ is a fixed compact subset of $B_{0}$. We show that $T$
maps $\left(A_{0},A_{1}\right)_{\theta,p}$ compactly into $\left(B_{0},B_{1}\right)_{\theta,p}$
for each $\theta\in(0,1)$ and $p\in[1,\infty]$.
\end{abstract}
This paper can be considered as a sequel to \cite{ci1}. Since the
introduction to \cite{ci1} provides exactly the background and references
which are relevant here, we shall mostly avoid repetition and request
the reader to be familiar with the first three pages of that paper.
However it will be convenient here to restate the following result
which appeared as Theorem 2 on p.~2 of \cite{ci1}.
\begin{thm}
\label{thm:simple}Let $\left(A_{0},A_{1}\right)$ and $\left(B_{0},B_{1}\right)$
be Banach couples. Suppose that $A_{0}\subset A_{1}$. Let $T$ be
a (possibly nonlinear) map of $A_{1}$ into $B_{1}$ which satisfies
the following two properties:\begin{equation}
T\left(A_{0}\right)\subset B_{0}\mbox{ and }\left\Vert T(a)\right\Vert _{B_{0}}\le C_{0}\left\Vert a\right\Vert _{A_{0}}\,\mbox{for each }a\in A_{0}\,,\label{eq:sone}\end{equation}
and\begin{equation}
\left\Vert Ta-Ta'\right\Vert _{B_{1}}\le C_{1}\left\Vert a-a'\right\Vert _{A_{1}}\,\mbox{for all }a,a'\in A_{1}\,.\label{eq:stwo}\end{equation}
where $C_{0}$ and $C_{1}$ are positive constants.

Then $T$ maps the space $\left(A_{0},A_{1}\right)_{\theta,p}$ boundedly
into $\left(B_{0},B_{1}\right)_{\theta,p}$ for each $\theta\in(0,1)$
and $p\in[1,\infty]$, and satisfies the estimate \begin{equation}
\left\Vert Ta\right\Vert _{\left(B_{0},B_{1}\right)_{\theta,p}}\le C_{0}^{1-\theta}C_{1}^{\theta}\left\Vert a\right\Vert _{\left(A_{0},A_{1}\right)_{\theta,p}}\,\mbox{for all }a\in\left(A_{0},A_{1}\right)_{\theta,p}\,.\label{eq:apqm}\end{equation}

\end{thm}
Two paragraphs after stating this theorem we posed a question regarding
interpolation of compactness properties of $T$ and later showed that
the answer to it is negative. But now, in the present paper, we will
obtain an affirmative answer to a variant of that question. Following
in the footsteps of Arne Persson \cite{arnep}, we will only consider
Banach couples $\left(B_{0},B_{1}\right)$ which satisfy a certain
approximation condition which will be recalled below. Then, if instead
of merely requiring $T$ to map $A_{0}$ into $B_{0}$ compactly (i.e.,
to map bounded subsets of $A_{0}$ to relatively compact subsets of
$B_{0}$), we require this map to be compact in a more quantitative
or {}``uniform'' sense, also to be specified below, this suffices
to ensure that $T$ maps $\left(A_{0},A_{1}\right)_{\theta,p}$ compactly
into $\left(B_{0},B_{1}\right)_{\theta,p}$ for each $\theta\in(0,1)$
and each $p\in[1,\infty]$.

Theorem \ref{thm:simple} is well known and is a simple special case
of various more elaborate results due to Jacques-Louis Lions \cite{lions},
Jaak Peetre \cite{PeetreCluj} and Luc Tartar \cite{tartar}. For
the reader's convenience we recall its very simple and short proof.
Given $a\in A_{1}=A_{0}+A_{1}$ and $t>0$ and $\varepsilon>0$, we
choose a decomposition $a=a_{0}+a_{1}$ such that \[
\left\Vert a_{0}\right\Vert _{A_{0}}+\frac{tC_{1}}{C_{0}}\left\Vert a_{1}\right\Vert _{A_{1}}\le(1+\varepsilon)K\left(\frac{tC_{1}}{C_{0}},a;A_{0},A_{1}\right).\]
Then $Ta=Ta-Ta_{0}+Ta_{0}$ and \begin{eqnarray*}
K(t,Ta;B_{0},B_{1}) & \le & \left\Vert Ta\right\Vert _{B_{0}}+t\left\Vert Ta-Ta_{0}\right\Vert _{B_{1}}\\
 & \le & C_{0}\left\Vert a_{0}\right\Vert _{A_{0}}+tC_{1}\left\Vert a-a_{0}\right\Vert _{A_{1}}\\
 & = & C_{0}\left(\left\Vert a_{0}\right\Vert _{A_{0}}+\frac{tC_{1}}{C_{0}}\left\Vert a_{1}\right\Vert \right)\\
 & \le & C_{0}(1+\varepsilon)K\left(\frac{tC_{1}}{C_{0}},a;A_{0},A_{1}\right)\end{eqnarray*}
and this immediately implies that \begin{equation}
K(t,Ta;B_{0},B_{1})\le C_{0}K\left(\frac{tC_{1}}{C_{0}},a;A_{0},A_{1}\right).\label{eq:pqm}\end{equation}
Exactly as in the {}``classical'' case of linear operators, the
estimate (\ref{eq:apqm}) and the inclusion $T\left(\left(A_{0},A_{1}\right)_{\theta,p}\right)\subset\left(B_{0},B_{1}\right)_{\theta,p}$
are obvious and trivial consequences of (\ref{eq:pqm}), especially
after we observe that (\ref{eq:pqm}) is equivalent to \[
t^{-\theta}K(t,Ta;B_{0},B_{1})\le C_{0}^{1-\theta}C_{1}^{\theta}\left(\frac{tC_{1}}{C_{0}}\right)^{-\theta}K\left(\frac{tC_{1}}{C_{0}},a;A_{0},A_{1}\right).\]
$\qed$

In the work of Persson, and also in many subsequent papers, for example,
\cite{CP} and \cite{Cobos-LPLemma}, it was shown that compactness
results for the special case where the {}``range'' couple $\left(B_{0},B_{1}\right)$
satisfies $B_{0}=B_{1}$ can be the first step towards obtaining corresponding
compactness results for more general couples $\left(B_{0},B_{1}\right)$.
The prototype of these results for the case $B_{0}=B_{1}$ is a classical
result of Jacques-Louis Lions and Jaak Peetre, often referred to as
the {}``Lions--Peetre Lemma''. It appears as Théorème 2.2 on p.~37
of \cite{LP}. It is remarked in \cite{LP} that similar ideas about
compactness appear in the work \cite{gag} of Emilio Gagliardo. Here
too, a variant of the Lions--Peetre Lemma will be an important ingredient
for us. We present it as the following lemma. Its proof is quite similar
to the one for linear operators in \cite{LP}. Another quite similar
result, for the case where $T:A_{0}\to B_{0}$ is also Lipschitz and
where $A_{0}$ is not necessarily contained in $A_{1}$, is due to
Fernando Cobos and appears as part (i) of Theorem 2.1 on p.~274 of
his paper \cite{Cobos-LPLemma}.
\begin{lem}
\label{lem:giraffe} Suppose that the Banach couples $\left(A_{0},A_{1}\right)$
and $\left(B_{0},B_{1}\right)$ and the map $T:A_{1}\to B_{1}$ satisfy
all the hypotheses of Theorem \ref{thm:simple}. Suppose, furthermore,
that $T$ maps every bounded subset of $A_{0}$ to a relatively compact
subset of $B_{0}$. Then, for each $\theta\in(0,1)$, and each $p\in[1,\infty]$,
$T$ maps every bounded subset of $\left(A_{0},A_{1}\right)_{\theta,p}$
to a relatively compact subset of $B_{0}+B_{1}$.

In particular, if $B_{0}=B_{1}$, then $T$ maps every bounded subset
of $\left(A_{0},A_{1}\right)_{\theta,p}$ to a relatively compact
subset of $B_{0}$.
\end{lem}
\textit{Proof.} Suppose that $\left\{ a_{n}\right\} _{n\in\mathbb{N}}$
is an arbitrary bounded sequence in $\left(A_{0},A_{1}\right)_{\theta,p}$.
Then it is also a bounded sequence in $\left(A_{0},A_{1}\right)_{\theta,\infty}$.
So, for each $m\in\mathbb{Z}$, and some constant $C_{2}$ we can
write $a_{n}=u_{n}(m)+v_{n}(m)$ where \[
\left\Vert u_{n}(m)\right\Vert _{A_{0}}+2^{m}\left\Vert v_{n}(m)\right\Vert _{A_{1}}\le C_{2}2^{\theta m}\,.\]
Then \begin{eqnarray*}
 &  & \left\Vert Tu_{n}(m)\right\Vert _{B_{0}}+2^{m}\left\Vert Ta_{n}-Tu_{n}(m)\right\Vert _{B_{1}}\\
 & \le & \max\left\{ C_{0},C_{1}\right\} \left(\left\Vert u_{n}(m)\right\Vert _{A_{0}}+2^{m}\left\Vert a_{n}-u_{n}(m)\right\Vert _{A_{1}}\right)\le C_{3}\cdot2^{\theta m}.\end{eqnarray*}
Thus we have \begin{equation}
\left\Vert Ta_{n}-Tu_{n}(m)\right\Vert _{B_{1}}\le C_{3}\cdot2^{-(1-\theta)m}.\label{eq:gqq}\end{equation}

We can suppose, by passing to subsequences and using Cantor diagonalization,
that, for each fixed $m\in\mathbb{Z}$, the sequence $\left\{ Tu_{n}(m)\right\} _{n\in\mathbb{N}}$
converges in $B_{0}$ norm, as $n$ tends to $\infty$, to some element
in $B_{0}$ depending on $m$. More relevantly for us, this means
that \begin{equation}
\left\Vert Tu_{n}(m)-Tu_{k}(m)\right\Vert _{B_{0}}\le\rho(m,n,k),\label{eq:myy}\end{equation}
where $\lim_{n,k\to\infty}\rho(m,n,k)=0$ for each fixed $m$. (The
rate of convergence here may depend on $m$.)

For each $n$ and $k$ in $\mathbb{N}$ we have \[
Ta_{n}-Ta_{k}=(Ta_{n}-Tu_{n}(m))+(Tu_{n}(m)-Tu_{k}(m))+(Tu_{k}(m)-Ta_{k}).\]
So\begin{equation}
\begin{array}{cl}
 & \left\Vert Ta_{n}-Ta_{k}\right\Vert _{B_{0}+B_{1}}\\
\le & \left\Vert Ta_{n}-Tu_{n}(m)\right\Vert _{B_{1}}+\left\Vert Tu_{n}(m)-Tu_{k}(m)\right\Vert _{B_{0}}+\left\Vert Tu_{k}(m)-Ta_{k}\right\Vert _{B_{1}}\\
\le & C_{3}\cdot2^{-(1-\theta)m}+\rho(m,n,k)+C_{3}\cdot2^{-(1-\theta)m}.\end{array}\label{eq:oomy}\end{equation}
Given any $\varepsilon>0$, there exists $m=m(\varepsilon)$ for which
$2\cdot C_{3}\cdot2^{-(1-\theta)m}<\varepsilon/2$. Then, for this
choice of $m$, there exists $N=N(\varepsilon)$ such that $\rho(m,n,k)<\varepsilon/2$
for all $n,k\ge N(\varepsilon)$. Thus we deduce from (\ref{eq:oomy})
that $\left\Vert Ta_{n}-Ta_{k}\right\Vert _{B_{0}+B_{1}}<\varepsilon$
whenever $n,k\ge N(\varepsilon)$ and the proof is complete. $\qed$

Let us recall the condition {}``(H)'' on Banach couples $\left(B_{0},B_{1}\right)$
which was introduced by Arne Persson in \cite{arnep}. The following
formulation of it looks very slightly different but is obviously equivalent.

(H) For each compact subset $K$ of $B_{0}$, there exists a positive
constant $c(K)$ and a set $\mathcal{D}(K)$ of linear operators $P:B_{0}+B_{1}\to B_{0}\cap B_{1}$
which satisfy \[
\Vert Pb\Vert_{B_{j}}\le c(K)\Vert b\Vert_{B_{j}}\,\mbox{ for each }b\in B_{j}\,\mbox{ and }j=0,1\,.\]
Furthermore, for each $\varepsilon>0$, there exists $P_{\varepsilon}\in\mathcal{D}(K)$
such that \[
\Vert P_{\varepsilon}b-b\Vert_{B_{0}}<\varepsilon\,\mbox{ for each }b\in K\,.\]

\begin{rem}
Many Banach couples which occur in applications satisfy condition
(H). For example, essentially the same arguments as on pp.~217--219
of \cite{arnep} show that the couple $\left(L^{p_{0}},L^{p_{1}}\right)$
on an arbitrary underlying measure space satisfies condition (H) for
all $p_{0}\in[1,\infty)$ and $p_{1}\in[1,\infty]$, and that if the
measure space is finite, then this also holds when $p_{0}=\infty$.
\end{rem}
Here is our main result. Its proof is a natural extension of the proof
of the theorem on p.~217 of \cite{arnep}.
\begin{thm}
Suppose that the Banach couples $\left(A_{0},A_{1}\right)$ and $\left(B_{0},B_{1}\right)$
and the map $T:A_{1}\to B_{1}$ satisfy all the hypotheses of Theorem
\ref{thm:simple}. Suppose, furthermore, that the couple $\left(B_{0},B_{1}\right)$
has Arne Persson's property (H) and that there exists some compact
subset $K$ of $B_{0}$ for which \begin{equation}
Ta\in\left\Vert a\right\Vert _{A_{0}}K\,\mbox{for each }a\in A_{0}\,.\label{eq:newc}\end{equation}
Then $T$ maps every bounded subset of $\left(A_{0},A_{1}\right)_{\theta,p}$
to a relatively compact subset of $\left(B_{0},B_{1}\right)_{\theta,p}$
for each $\theta\in(0,1)$ and each $p\in[1,\infty]$.
\end{thm}
\textit{Proof.} The condition (\ref{eq:newc}) ensures that every
non zero element $a\in A_{0}$ satisfies $\frac{1}{\|a\|_{A_{0}}}Ta\in K$.
So, for each $\varepsilon>0$, Persson's property provides us with
a linear operator $P_{\varepsilon}\in\mathcal{D}(K)$ such that \begin{equation}
\left\Vert P_{\varepsilon}\left(\frac{1}{\|a\|_{A_{0}}}Ta\right)-\frac{1}{\|a\|_{A_{0}}}Ta\right\Vert _{B_{0}}<\varepsilon\,.\label{eq:brahms}\end{equation}
We shall show that the map $P_{\varepsilon}T-T$ satisfies all the
conditions of Theorem \ref{thm:simple} with $C_{0}$ and $C_{1}$
replaced by certain other constants:

Since $P_{\varepsilon}$ is linear, the condition (\ref{eq:brahms})
implies that $\|\left(P_{\varepsilon}T-T\right)a\|_{B_{0}}\leq\varepsilon\|a\|_{A_{0}}$
holds for each non zero $a\in A_{0}$. This estimate also holds for
$a=0$, since $T0=0$ (by (\ref{eq:newc}) or (\ref{eq:sone})).

Each pair of elements $a,a'\in A_{1}$ satisfies \begin{eqnarray*}
\|\left(P_{\varepsilon}T-T\right)a'-\left(P_{\varepsilon}T-T\right)a\|_{B_{1}} & \le & \|P_{\varepsilon}\left(Ta'-Ta\right)\|_{B_{1}}+\|Ta'-Ta\|_{B_{1}}\\
 & \le & \left(c\left(K\right)+1\right)\|Ta'-Ta\|_{B_{1}}\\
 & \le & C_{1}\left(c\left(K\right)+1\right)\|a'-a\|_{A_{1}}\,.\end{eqnarray*}

The preceding estimates enable us to apply Theorem \ref{thm:simple}
with $C_{0}$ and $C_{1}$ replaced by $\varepsilon$ and $C_{1}(c(K)+1)$
to obtain that \begin{equation}
\|\left(P_{\varepsilon}T-T\right)a\|_{\left(B_{0},B_{1}\right)_{\theta,p}}\leq\varepsilon^{1-\theta}\left(C_{1}\left(c\left(K\right)+1\right)\right)^{\theta}\|a\|_{\left(A_{0},A_{1}\right)_{\theta,p}}\,\mbox{ for all }a\in\left(A_{0},A_{1}\right)_{\theta,p}\,.\label{eq:bach}\end{equation}

Our next step is to apply Lemma \ref{lem:giraffe} to the map $P_{\varepsilon}T$,
but with the couple $\left(B_{0},B_{1}\right)$ replaced by $\left(\left(B_{0},B_{1}\right)_{\theta,p},\left(B_{0},B_{1}\right)_{\theta,p}\right)$.
As pointed out by Persson, the Closed Graph Theorem ensures that $P_{\varepsilon}$
is a bounded map from $B_{j}$ to $B_{0}\cap B_{1}$ for $j=0,1$.
Therefore $P_{\varepsilon}:B_{j}\to\left(B_{0},B_{1}\right)_{\theta,p}$
is also bounded and we have $P_{\varepsilon}T\left(A_{j}\right)\subset\left(B_{0},B_{1}\right)_{\theta,p}$.

Each element $a\in A_{0}$ satisfies \[
\|P_{\varepsilon}Ta\|_{\left(B_{0},B_{1}\right)_{\theta,p}}\leq\left\Vert P_{\varepsilon}\right\Vert _{B_{0}\to\left(B_{0},B_{1}\right)_{\theta,p}}\|Ta\|_{B_{0}}\leq C_{0}\left\Vert P_{\varepsilon}\right\Vert _{B_{0}\to\left(B_{0},B_{1}\right)_{\theta,p}}\|a\|_{A_{0}}\]
 and each pair of elements $a,a'\in A_{1}$ satisfies \begin{eqnarray*}
 &  & \|P_{\varepsilon}Ta'-P_{\varepsilon}Ta\|_{\left(B_{0},B_{1}\right)_{\theta,p}}\\
 & \le & \left\Vert P_{\varepsilon}\right\Vert _{B_{1}\to\left(B_{0},B_{1}\right)_{\theta,p}}\|Ta'-Ta\|_{B_{1}}\leq C_{1}\left\Vert P_{\varepsilon}\right\Vert _{B_{1}\to\left(B_{0},B_{1}\right)_{\theta,p}}\|a'-a\|_{A_{1}}\,.\end{eqnarray*}
Let $M$ be a bounded subset of $A_{0}$. $T$ maps $M$ to a relatively
compact subset of $B_{0}$, and thus $P_{\varepsilon}T$ maps $M$
to a relatively compact subset of $\left(B_{0},B_{1}\right)_{\theta,p}$
since $P_{\varepsilon}:B_{0}\to\left(B_{0},B_{1}\right)_{\theta,p}$
is continuous.

In view of all these properties of $P_{\varepsilon}T$, we can apply
the last part of the statement of Lemma \ref{lem:giraffe} to obtain
that $P_{\varepsilon}T$ maps each bounded subset of $\left(A_{0},A_{1}\right)_{\theta,p}$
to a relatively compact subset of $\left(B_{0},B_{1}\right)_{\theta,p}$.

For the final step of the proof we let $M$ be an arbitrary bounded
subset of $\left(A_{0},A_{1}\right)_{\theta,p}$. We have to show
that $T(M)$ is relatively compact in $\left(B_{0},B_{1}\right)_{\theta,p}$.
In fact we will show that it is totally bounded in $\left(B_{0},B_{1}\right)_{\theta,p}$.

Let $\varepsilon$ be an arbitrary positive number and choose $\varepsilon_{0}>0$
such that \[
\varepsilon_{0}^{1-\theta}\left(C_{1}\left(C\left(K\right)+1\right)\right)^{\theta}\|a\|_{\left(A_{0},A_{1}\right)_{\theta,p}}<\frac{\varepsilon}{2}\;\mbox{ for all }a\in M\,.\]
In view of (\ref{eq:bach}), this ensures that \begin{equation}
\|\left(P_{\varepsilon_{0}}T-T\right)a\|_{\left(B_{0},B_{1}\right)_{\theta,p}}<\frac{\varepsilon}{2}\,\mbox{ for all }a\in M\,.\label{eq:mozart}\end{equation}
Since $P_{\varepsilon_{0}}T\left(M\right)$ is relatively compact
and thus totally bounded in $\left(B_{0},B_{1}\right)_{\theta,p}$,
the inclusion $P_{\varepsilon_{0}}T\left(M\right)\subset{\displaystyle \bigcup_{i=1}^{n}}B\left(x_{i},\frac{\varepsilon}{2}\right)$
holds for some finite subset $\left\{ x_{i}\right\} _{i=1}^{n}$ of
$\left(B_{0},B_{1}\right)_{\theta,p}$. (Here of course $B(x,r)$
denotes the open ball in $\left(B_{0},B_{1}\right)_{\theta,p}$ of
radius $r$ centred at $x$.) This, together with the inclusion $T\left(M\right)\subset\left(T-P_{\varepsilon_{0}}T\right)\left(M\right)+P_{\varepsilon_{0}}T\left(M\right)$
and (\ref{eq:mozart}), gives us that $T\left(M\right)\subset{\displaystyle \bigcup_{i=1}^{n}}B\left(x_{i},\varepsilon\right)$
and completes the proof. $\qed$

\end{document}